\documentclass[conference]{IEEEtran2015}
\IEEEoverridecommandlockouts
\usepackage{cite}
\usepackage{amsmath,mathtools}
\usepackage{amssymb}
\usepackage[usenames,dvipsnames]{pstricks}
\usepackage{pstricks-add}
\usepackage{pst-grad}
\usepackage{pst-plot} 
\usepackage{pst-eucl} 
\psset{algebraic}
\usepackage{pifont}
\usepackage{theorem} 
\usepackage{euscript}
\usepackage{exscale,relsize}
\usepackage{epic,eepic}
\usepackage{tikz}
\usepackage{pgfplots}
\usepackage{makeidx}
\usepackage{multicol}
\usepackage{srcltx}         %for forward references in xdvi
\usepackage{xcolor}
\usepackage{url}
\usepackage{charter}
\usepackage{mathrsfs} %pour mathscr
\usepackage{graphicx}
\definecolor{dred}{HTML}{D90404}
\definecolor{orng}{HTML}{D35400}

\definecolor{refkey}{rgb}{0,0.6,0.0}
\definecolor{Brown}{rgb}{0.45,0.0,0.05}
\definecolor{dgreen}{rgb}{0.00,0.49,0.00}
\definecolor{dblue}{rgb}{0,0.08,0.75}
\definecolor{lblue}{rgb}{0,0.7,0.75}
\renewcommand{\leq}{\ensuremath{\leqslant}}

\newcommand{\scal}[2]{{\langle{{#1}\mid{#2}}\rangle}}

\newcommand{\menge}[2]{\big\{{#1}~\big |~{#2}\big\}}

\newcommand{\HH}{\ensuremath{{\mathcal H}}}
\newcommand{\GG}{\ensuremath{{\mathcal G}}}

\newcommand{\emp}{\ensuremath{{\varnothing}}}

\newcommand{\Id}{\ensuremath{\operatorname{Id}}\,}

\newcommand{\RR}{\ensuremath{\mathbb{R}}}

\newcommand{\soft}[1]{\ensuremath{{\operatorname{soft}}_{{#1}}\,}}
\newcommand{\hard}[1]{\ensuremath{{\operatorname{hard}}_{{#1}}\,}}

\newcommand{\NN}{\ensuremath{\mathbb N}}

\newcommand{\pinf}{\ensuremath{{+\infty}}}

\newcommand{\proj}{\ensuremath{\text{\rm proj}}}

\newcommand{\Fix}{\ensuremath{\text{\rm Fix}\,}}

\newcommand{\dft}{\ensuremath{\operatorname{DFT}}\,}

\newcommand{\sign}{\ensuremath{\text{\rm sign}}}

\newcommand{\rzeroun}{\ensuremath{\left]0,1\right]}}   
   
\newtheorem{theorem}{Theorem}
\theoremstyle{plain}{\theorembodyfont{\rmfamily}
}

\theoremstyle{plain}{\theorembodyfont{\rmfamily}
\newtheorem{assumption}[theorem]{Assumption}}
\theoremstyle{plain}{\theorembodyfont{\rmfamily}
}
\usepackage{enumitem}
\renewcommand\theenumi{(\roman{enumi})}
\renewcommand\theenumii{\alph{enumii})}

\begin{document}
\title{A FIXED POINT FRAMEWORK FOR RECOVERING SIGNALS FROM NONLINEAR 
TRANSFORMATIONS\thanks{The work of P. L. Combettes was supported by 
the National Science Foundation under grant CCF-1715671 and the
work of Z. C. Woodstock was supported by the National Science 
Foundation under grant DGE-1746939.}}
\author{\IEEEauthorblockN{%1\textsuperscript{st} 
Patrick L. Combettes}
\IEEEauthorblockA{\textit{Department of Mathematics}\\
\textit{North Carolina State University}\\
Raleigh, NC 27695-8205, USA}
\and
\IEEEauthorblockN{%2\textsuperscript{nd} 
Zev C. Woodstock}
\IEEEauthorblockA{\textit{Department of Mathematics}\\
\textit{North Carolina State University}\\
Raleigh, NC 27695-8205, USA}}

\maketitle

\begin{abstract}
We consider the problem of recovering a signal from
nonlinear transformations, under convex constraints modeling
\emph{a priori} information. Standard feasibility and optimization
methods are ill-suited to tackle this problem due to the
nonlinearities. We show that, in many common
applications, the transformation model can be
associated with fixed point equations involving firmly
nonexpansive operators. In turn, the
recovery problem is reduced to a tractable common fixed point
formulation, which is solved efficiently by a provably convergent,
block-iterative algorithm. Applications to signal and image
recovery are demonstrated. Inconsistent problems are also
addressed.
\end{abstract}


\section{Introduction}
\label{sec:1}
Under consideration is the general problem of recovering an
original signal $\overline{x}$ in a Euclidean space
$\HH$ from a finite number of transformations $(r_k)_{k\in K}$
of the form
\begin{equation}
\label{e:1}
r_k=R_k\overline{x}\in\GG_k,
\end{equation}
where $R_k\colon\HH\to\GG_k$ is an operator mapping the solution
space $\HH$ to the Euclidean space $\GG_k$. 
In addition to these transformations, some \emph{a priori} 
constraints on $\overline{x}$ are available in the
form of a finite family of closed convex subsets $(C_j)_{j\in J}$
of $\HH$ \cite{Youl82,Aiep96,Rzep18,Tofi16,Trus84}. Altogether, the 
recovery problem is to
\begin{equation}
\label{e:2}
\text{find}\;\;x\in\bigcap_{j\in J}C_j\;\;\text{such that}\;\;
(\forall k\in K)\quad R_kx=r_k.
\end{equation}
One of the most classical instances of this formulation was
proposed by Youla in \cite{Youl78}, namely 
\begin{equation}
\label{e:3}
\text{find}\;\;x\in V_1\;\;\text{such that}\;\;\proj_{V_2}x=r_2,
\end{equation}
where $V_1$ and $V_2$ are vector subspaces of $\HH$ and
$\proj_{V_2}$ is the projection operator onto $V_2$. As shown in 
\cite{Youl78}, 
\eqref{e:3} covers many basic signal processing problems, such as
band-limited extrapolation or image reconstruction from diffraction
data, and it can be solved with a simple alternating projection
algorithm. The extension of \eqref{e:3} to recovery problems with
several transformations modeled as linear projections 
$r_k=\proj_{V_k}\overline{x}$ is discussed in
\cite{Joat10,Reye13}.
More broadly, if the operators $(R_k)_{k\in K}$
are linear, reliable algorithms are available to solve \eqref{e:2}.
In particular, since the associated constraint set is 
an affine subspace with an explicit projection operator, standard 
feasibility
algorithms can be used \cite{Aiep96}. Alternatively, proximal
splitting methods can be considered; see \cite{MaPr18} and its
references.

In the present paper we consider the general situation in
which the operators $(R_k)_{k\in K}$ in \eqref{e:1} are not 
necessarily linear, a stark departure from common assumptions in 
signal recovery problems. Examples of such nonlinearly generated
data $(r_k)_{k\in K}$ in \eqref{e:1} include
hard-thresholded wavelet coefficients of $\overline{x}$, the 
positive part of the Fourier transform of $\overline{x}$, a
mixture of best approximations of $\overline{x}$ from 
closed convex sets, a maximum a posteriori denoised version of
$\overline{x}$, or measurements of $\overline{x}$ acquired through
nonlinear sensors. 

A significant difficulty one faces in the nonlinear context 
is that the constraint \eqref{e:1} is typically not representable
by an exploitable convex constraint; see, e.g., 
\cite{Blum13,Cast19}. As a result, finding a solution to
\eqref{e:2} with a provenly convergent and numerically efficient
algorithm is a challenging task. In particular, standard convex 
feasibility algorithms are not applicable. Furthermore, 
variational relaxations involving a penalty of the type 
$\sum_{k\in K}\phi_k(\|R_kx-r_k\|)$ typically lead to nonconvex
problems, even for choices as basic as $\phi_k=|\cdot|^2$ and $R_k$
taken as the projection operator onto a closed convex set. 

Our strategy to solve \eqref{e:2} is to forego the feasibility and 
optimization approaches in favor of the flexible and unifying 
framework of fixed point theory. Our first contribution is
to show that, while $R_k$ in \eqref{e:1} may be a very badly 
conditioned (possibly discontinuous) operator, common transformation 
models can be reformulated as fixed point equations
with respect to an operator with much better properties,
namely a firmly nonexpansive operator. Next, using a
suitable modeling of the constraint sets $(C_j)_{j\in J}$, we 
rephrase \eqref{e:2} as an equivalent common fixed point
problem and solve it with a reliable and efficient extrapolated
block-iterative fixed point algorithm. This strategy is outlined
in Section~\ref{sec:2}, where we also provide the algorithm. In
Section~\ref{sec:3}, we present several numerical illustrations of
the proposed framework to nonlinear signal and image recovery.
Finally, inconsistent problems are addressed in Section~\ref{sec:4}.

\section{Fixed point model and algorithm}
\label{sec:2}
For background on the tools from fixed point theory and convex
analysis used in this section, we refer the reader to
\cite{Livre1}. Let us first recall that an operator 
$T\colon\HH\to\HH$ is firmly nonexpansive if
\begin{multline}
\label{e:10}
(\forall x\in\HH)(\forall y\in\HH)\quad\|Tx-Ty\|^2\leq\\
\|x-y\|^2-\|(\Id-T)x-(\Id-T)y\|^2,
\end{multline}
and firmly quasinonexpansive if
\begin{equation}
\label{e:11}
(\forall x\in\HH)(\forall y\in\Fix T)\quad
\scal{y-Tx}{x-Tx}\leq 0,
\end{equation}
where $\Fix T=\menge{x\in\HH}{Tx=x}$. Finally, the subdifferential 
of a convex function $f\colon\HH\to\RR$ at $x\in\HH$ is
\begin{equation}
\partial f(x)\!=\!\menge{u\in\HH\!}{\!(\forall y\in\HH)
\,\scal{y-x}{u}
+f(x)\!\leq\!f(y)}.
\end{equation}

As discussed in Section~\ref{sec:1}, the transformation model
\eqref{e:1} is too general to make finding a solution to 
\eqref{e:2} via a provenly convergent method possible. We
therefore assume the following.

\begin{assumption}
\label{a:1}
The problem \eqref{e:2} has at least one solution, $J\cap K=\emp$, 
and the following hold:
\begin{enumerate}
\item
\label{a:1iii}
For every $k\in K$, $S_k\colon\GG_k\to\HH$ is an operator such 
that $S_k\circ R_k$ is firmly nonexpansive and
\begin{equation}
\label{e:ai}
\bigg(\forall x\in\bigcap_{j\in J}C_j\bigg)\;
S_k(R_kx)=S_kr_k\;\Rightarrow\; R_kx=r_k.
\end{equation}
\item
\label{a:1i}
For every $j\in J_1\subset J$, the operator $\proj_{C_j}$ is
easily implementable.
\item
\label{a:1ii}
For every $j\in J\smallsetminus J_1$, 
$f_j\colon\HH\to\RR$ is a convex function such that 
$C_j=\menge{x\in\HH}{f_j(x)\leq 0}$.
\end{enumerate}
\end{assumption}

In view of Assumption~\ref{a:1}\ref{a:1iii}, let us replace
\eqref{e:2} by the equivalent problem 
\begin{equation}
\label{e:22}
\text{find}\;x\in\bigcap_{j\in J}C_j\;\text{such that}\;
(\forall k\in K)\;S_k(R_kx)=S_kr_k.
\end{equation}
Concrete examples of suitable operators 
$(S_k)_{k\in K}$ will be given in Section~\ref{sec:3} 
(see also \cite{Ibap20}). The motivation behind 
\eqref{e:22} is that it leads to a tractable fixed point
formulation. To see this, set
\begin{equation}
\label{e:21}
(\forall k\in K)\quad T_k=S_kr_k+\Id-S_k\circ R_k
\end{equation}
and let $x\in\bigcap_{j\in J}C_j$.
Then, for every $k\in K$, \eqref{e:1} 
$\Leftrightarrow$ $S_k(R_k{x})=S_kr_k$ $\Leftrightarrow$ 
${x}=S_kr_k+{x}-S_k(R_k{x})$
$\Leftrightarrow$ ${x}\in\Fix T_k$. A key observation at
this point is that \eqref{e:10} implies that the operators 
$(T_k)_{k\in K}$ are firmly nonexpansive, hence firmly 
quasinonexpansive. 

If $j\in J_1$, per Assumption~\ref{a:1}\ref{a:1i}, the set 
$C_j$ will be activated in the algorithm through the use of the
operator $T_j=\proj_{C_j}$, which is firmly nonexpansive 
\cite[Proposition~4.16]{Livre1}. 
On the other hand, if $j\in J\smallsetminus J_1$, 
the convex inequality representation of
Assumption~\ref{a:1}\ref{a:1ii} will lead
to an activation of $C_j$ through its subgradient projector. 
Recall that the subgradient projection of $x\in\HH$ onto $C_j$
relative to $u_j\in\partial f_j(x)$ is 
\begin{equation}
\label{e:7}
T_jx=
\begin{cases}
x-\dfrac{f_j(x)}{\|u_j\|^2}u_j,&\text{if}\;\;f_j(x)>0;\\
x,&\text{if}\;\;f_j(x)\leq 0,
\end{cases}
\end{equation}
and that $T_j$ is firmly quasinonexpansive, with $\Fix T_j=C_j$
\cite[Proposition~29.41]{Livre1}.
The advantage of the subgradient projector onto $C_j$ is that, 
unlike the exact projector, it does not require solving a
nonlinear best approximation problem, which makes it much easier
to implement in the presence of convex inequality constraints
\cite{Imag97}. Altogether, \eqref{e:2} is equivalent to the 
common fixed point problem 
\begin{equation}
\label{e:23}
\text{find}\;\;x\in\bigcap_{i\in J\cup K}\Fix T_i,
\end{equation}
where each $T_i$ is firmly quasinonexpansive. This allows us to
solve \eqref{e:2} as follows.

\begin{theorem}{\rm\cite{Ibap20}}
\label{t:1}
Consider the setting of problem \eqref{e:2} under
Assumption~\ref{a:1}. Let
$x_0\in\HH$, let $0<\varepsilon<1/\text{\rm card}(J\cup K)$, and
set $(\forall k\in K)$ $p_k=S_kr_k$ and $F_k=S_k\circ R_k$. Iterate
\begin{equation}
\label{e:alg}
\hskip -0.6mm
\begin{array}{l}
\text{for}\;\;n=0,1,\ldots\\
\left\lfloor
\begin{array}{l}
\emp\neq I_n\subset J\cup K\\
\{\omega_{i,n}\}_{i\in I_n}\subset[\varepsilon,1],\;
\sum_{i\in I_n}\omega_{i,n}=1\\
\text{for every}\;\;i\in I_n\\
\left\lfloor
\begin{array}{l}
\text{if}\;\;i\in J_1\\
\left\lfloor
\begin{array}{l}
y_{i,n}=\proj_{C_i}x_n-x_n\\
\end{array}
\right.\\
\text{if}\;\;i\in J\smallsetminus J_1\\
\left\lfloor
\begin{array}{l}
u_{i,n}\in\partial f_i(x_n)\\
y_{i,n}=
\begin{cases}
-\dfrac{f_i(x_n)}{\|u_{i,n}\|^2}u_{i,n}
&\text{if}\;f_i(x_n)>0\\
0,&\text{if}\;f_i(x_n)\leq 0
\end{cases}
\end{array}
\right.\\
\text{else}\\
\left\lfloor
y_{i,n}=p_i-F_ix_n
\right.\\
\nu_{i,n}=\|y_{i,n}\|\\
\end{array}
\right.\\
\nu_n=\sum_{i\in I_n} \omega_{i,n} \nu_{i,n}^2\\
\text{if}\;\nu_n=0\\
\left\lfloor
\begin{array}{l}
x_{n+1}=x_n
\end{array}
\right.\\
\text{else}\\
\left\lfloor
\begin{array}{l}
y_n=\sum_{i\in I_n}\omega_{i,n} y_{i,n}\\
\Lambda_n=\nu_n/\|y_n\|^2\\
\lambda_n\in[\varepsilon,(2-\varepsilon)\Lambda_n]\\
x_{n+1}=x_n+\lambda_n y_n.
\end{array}
\right.\\[6.8mm]
\end{array}
\right.\\
\end{array}
\end{equation}
Suppose that there exists an integer $M>0$ such 
that
\begin{equation}
\label{e:33}
(\forall n\in\NN)\quad\bigcup_{m=0}^{M-1}I_{n+m}=J\cup K. 
\end{equation}
Then $(x_n)_{n\in\NN}$ converges to a solution to \eqref{e:2}.
\end{theorem}

When $K=\emp$, \eqref{e:alg} coincides with the 
extrapolated method of parallel subgradient projections (EMOPSP)
of \cite{Imag97}. It has in addition the ability to
incorporate the constraints \eqref{e:1}, while maintaining the
attractive features of EMOPSP. First, it can process
the operators in blocks of variable size. 
The control scheme \eqref{e:33} just 
imposes that every operator be activated at least once within any
$M$ consecutive iterations. Second, because the extrapolation
parameters $(\Lambda_n)_{n\in\NN}$ can attain
large values in $\left[1,\pinf\right[$, large steps are possible, 
which lead to fast convergence compared to standard relaxation
schemes, where $\Lambda_n\equiv 1$.

\section{Applications}
\label{sec:3}
We illustrate several instances of \eqref{e:2}, develop tractable
reformulations of the form \eqref{e:22}, and solve them using
\eqref{e:alg}, where $x_0=0$ and the relaxation strategy
is that recommended in \cite[Chapter~5]{Aiep96}, namely
\begin{equation}
\label{e:lsteps}
(\forall n\in\NN)\quad\lambda_n=
\begin{cases}
\Lambda_n/2,&\text{if}\;\;n=0\mod 3;\\
1.99\Lambda_n,&\text{otherwise.}
\end{cases}
\end{equation}

\subsection{Restoration from distorted signals}
\label{sec:31}
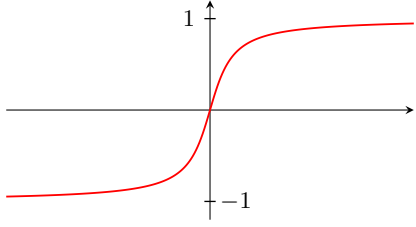
\begin{figure}[t]
\centering
\begin{tabular}{@{}c@{}c@{}}
\begin{tikzpicture}
\begin{axis}[
  width=7cm,
  height=4.5cm,
  ymin=-1.2,ymax=1.2,
  axis lines=middle,
  clip=false,
  xtick style={draw=none},
  xticklabels=\empty,
  ytick={1},
  yticklabels={$1$},
  ticklabel style={font=\small},
  tick style={color=black, line width=0.05em},
  extra y ticks={-1},
  extra y tick style={y tick label style={right,
  xshift=0.25em},font=\small},
  extra y tick labels={$-1$},
]
%Distortion
\addplot [line width=0.7pt,mark=none,red,
samples=250,domain=-1.2:1.2] {2*rad(atan(10*x))/pi};
\end{axis}
\end{tikzpicture}\\
\end{tabular} 
\caption{Distortion operator $\theta_{3}$ in Section~\ref{sec:31}.}
\label{fig:2'}
\end{figure}

The goal is to recover the original form of the $N$-point
($N=2048$) signal $\overline{x}$ from the following (see
Fig.~\ref{fig:1}):
\begin{itemize}
\item 
A bound $\gamma_1$ on the energy of the finite differences of
$\overline{x}$, namely $\|D\overline{x}\|\leq\gamma_1$, where 
$D\colon(\xi_i)_{i\in\{0,\ldots,N-1\}}\mapsto
(\xi_{i+1}-\xi_i)_{i\in\{0,\ldots,N-2\}}$. The bound is given 
from prior information as $\gamma_1=1.17$. 
\item
A distortion $r_2=R_2\overline{x}$,
where $R_2$ clips componentwise to $[-\gamma_2,\gamma_2]$ 
($\gamma_2=0.1$) \cite[Section~10.5]{Tarr18}.
\item
A distortion $r_3=R_3\overline{x}$ of a low-pass
version of $\overline{x}$, where
$R_3=Q_{3}\circ L_{3}$. Here $L_{3}$ bandlimits by zeroing all
but the $83$ lowest-frequency coefficients of the Discrete
Fourier Transform, and $Q_{3}$ induces componentwise distortion
via the operator \cite[Section~10.6]{Tarr18}
$\theta_3=(2/\pi)\arctan(\gamma_{3}\:\cdot)$,
where $\gamma_{3}=10$ (see Fig.~\ref{fig:2'}).
\end{itemize}
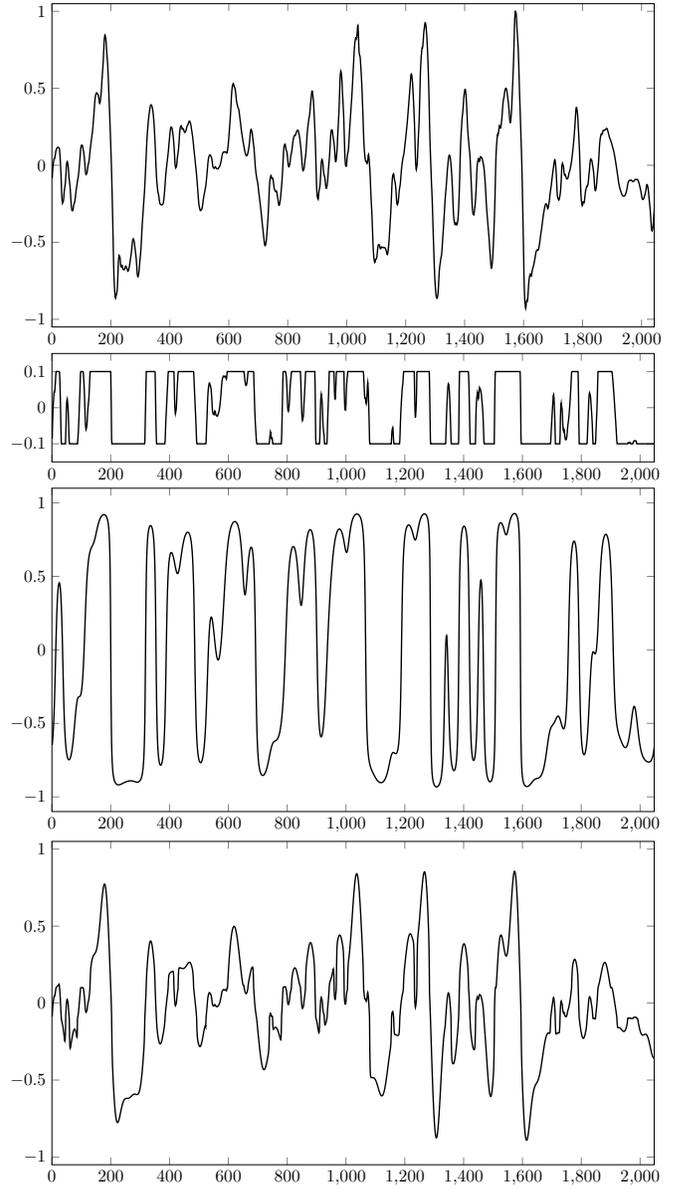
\begin{figure}[t]
\begin{tikzpicture}[scale=0.65]
\begin{axis}[height=8.2cm,width=13.9cm, legend cell align={left},
xmin =0, xmax=2043, ymin=-1.05, ymax=1.05]
\addplot [thick, mark=none, color=black] table[x={n},
y={xtrue}] {figures/ex1/ex1.txt};
\end{axis}
\end{tikzpicture}\\
\begin{tikzpicture}[scale=0.65]
\begin{axis}[height=3.8cm,width=13.9cm, legend cell align={left},
xmin=0, xmax=2048, ymin=-.15, ymax=.15]
\addplot [thick, mark=none, color=black] table[x={n},
y={sat}] {figures/ex1/ex1.txt};
\end{axis}
\end{tikzpicture}\\
\begin{tikzpicture}[scale=0.65]
%\hskip 1.75mm
\begin{axis}[height=8.2cm,width=13.9cm, legend cell align={left},
xmin =0, xmax=2048, ymin=-1.1, ymax=1.1]
\addplot [thick, mark=none, color=black] table[x={n},
y={distBL}] {figures/ex1/ex1.txt};
\end{axis}
\end{tikzpicture}\\
\begin{tikzpicture}[scale=0.65]
\begin{axis}[height=8.2cm,width=13.9cm, legend cell align={left},
xmin =0, xmax=2048, ymin=-1.05, ymax=1.05]
\addplot [thick, mark=none, color=black] table[x={n},
y={xnew}] {figures/ex1/ex1.txt};
\end{axis}
\end{tikzpicture}
\caption{Signals in Section~\ref{sec:31}. Top to bottom:
original signal $\overline{x}$, 
distorted signal $r_2$,
distorted signal $r_3$,
recovered signal.}
\label{fig:1}
\end{figure} 
The solution space is the standard Euclidean space
$\HH=\mathbb{R}^N$. To formulate the recovery problem as an 
instance of \eqref{e:2}, set $J=\{1\}$, $J_1=\emp$, $K=\{2,3\}$,
and $C_1=\menge{x\in\HH}{f_1(x)\leq 0}$, where
$f_1=\|D\cdot\|-\gamma_1$. Then the objective is to
\begin{multline}
\label{e:1d}
\text{find}\;\;x\in C_1\;\;\text{such that}\;\;R_2x=r_2\;\;
\text{and}\;\;R_3x=r_3.
\end{multline}
Next, let us verify that Assumption~\ref{a:1}\ref{a:1iii} is
satisfied. On the one hand, since $R_2$ is the projection onto
the closed convex set $[-\gamma_2,\gamma_2]^{N}$, it is firmly
nonexpansive, so we set
$S_2=\Id$. On the other hand, if we set
$S_3=\gamma_{3}^{-1}L_3$, then $S_3\circ R_3$ is
firmly nonexpansive and satisfies \eqref{e:ai} \cite{Ibap20}. 
We thus obtain an instance of \eqref{e:22}, to which we
apply \eqref{e:alg} with \eqref{e:lsteps} and
$(\forall n\in\mathbb{N})$
$I_n=J\cup K$ and $(\forall i\in I_n)$ $\omega_{i,n}=1/3$.
The recovered signal shown in Fig.~\ref{fig:1} effectively
incorporates the information from the prior constraint and 
the nonlinear distortions.

\subsection{Reconstruction from thresholded scalar products}
\label{sec:32}
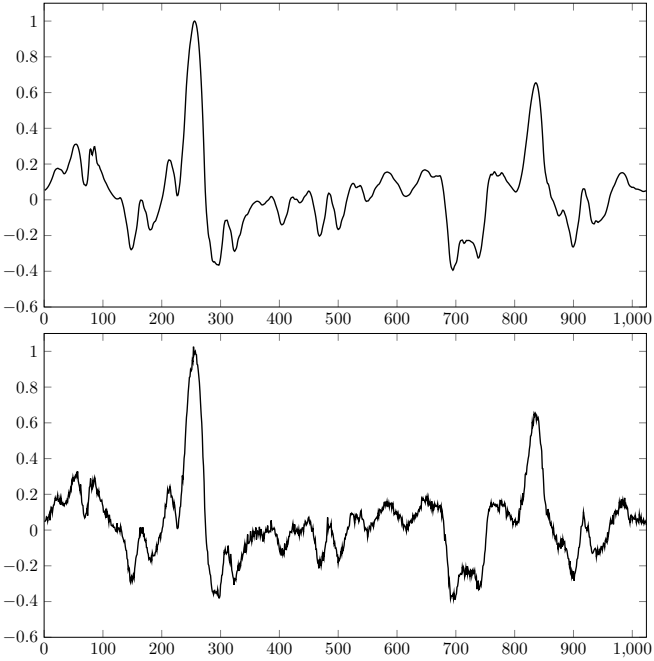
\begin{figure}[t]
\begin{tikzpicture}[scale=0.65]
\begin{axis}[height=7.8cm,width=13.9cm, legend cell align={left},
xmin =0, xmax=1024, ymin=-0.6, ymax=1.1]
\addplot [thick, mark=none, color=black] table[x={n},
y={xtrue}] {figures/ex2/ex2.txt};
\end{axis}
\end{tikzpicture}
\begin{tikzpicture}[scale=0.65]
\begin{axis}[height=7.8cm,width=13.9cm, legend cell align={left},
xmin =0, xmax=1024, ymin=-0.6, ymax=1.1]
\addplot [thick, mark=none, color=black] table[x={n}, y={xnew}]
{figures/ex2/ex2.txt};
\end{axis}
\end{tikzpicture}
\caption{Original signal $\overline{x}$ (top) and recovery 
(bottom) in Section~\ref{sec:32}.}
\label{fig:sig}
\end{figure}

\begin{figure}[b]
\centering
\begin{tikzpicture}
\begin{axis}[
width=7cm,
height=4.5cm,
xmin=-2.5,
xmax=2.5,
axis lines=middle,
clip=false,
xtick={0.7},
ytick style={draw=none},
xticklabels={$\gamma$},
yticklabels=\empty,
ticklabel style={font=\small},
extra x ticks={-0.7},
extra x tick style={x tick label style={above,
yshift=0.2em,xshift=-.2em},font=\small},
extra x tick labels={$-\gamma$},
tick style={color=black, line width=0.07em},
]
%Theta thresholder
\addplot [line width=0.7pt,mark=none,red,
samples=250,domain=0.7:2.5] {sqrt(x^2-0.7^2)};
\addplot [line width=0.7pt,mark=none,red,
samples=250,domain=-2.5:-0.7] {-sqrt(x^2-0.7^2)};
\addplot [line width=0.7pt,mark=none,red,
samples=250,domain=-0.7:0.7] {0.02};
\addplot [line width=0.7pt,mark=none,blue,
samples=250,domain=0.7:2.5] {x-0.7};
\addplot [line width=0.7pt,mark=none,blue,
samples=250,domain=-2.5:-0.7] {x+0.7};
\addplot [line width=0.7pt,mark=none,blue,
samples=250,domain=-0.7:0.7] {-0.02};
\end{axis}
\end{tikzpicture}
\caption{The thresholder \eqref{e:22-1} of \cite{Taov00} (red) and
the soft thresholder (blue) used in Section~\ref{sec:32}.}
\label{fig:thr}
\end{figure}
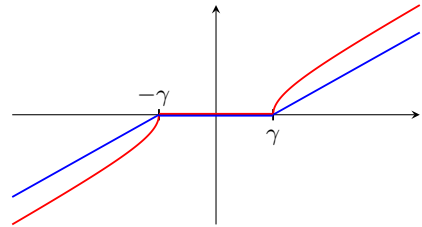

The goal is to recover the original form of the $N$-point
($N=1024$) signal $\overline{x}$ shown in Fig.~\ref{fig:sig} 
from thresholded scalar products $(r_k)_{k\in K}$ given by
\begin{multline}
\label{e:94}
(\forall k\in K)\quad r_k=R_k\overline{x},\quad\text{with}\\
R_k\colon\HH\to\RR\colon x\mapsto Q_{\gamma}\scal{x}{e_k},
\end{multline}
where 
\begin{itemize}
\item
$(e_k)_{k\in K}$ is a collection of normalized vectors in $\RR^N$
with zero-mean i.i.d. entries. 
\item
$Q_{\gamma}$ ($\gamma=0.05$) is the thresholding operator
\begin{equation}
\label{e:22-1}
Q_{\gamma}\colon\xi\mapsto
\begin{cases}
\sign(\xi)\sqrt{\xi^2-\gamma^2},
&\text{if}\;\;|\xi|>\gamma;\\
0,&\text{if}\;\;|\xi|\leq\gamma
\end{cases} 
\end{equation}
of \cite{Taov00} (see Fig.~\ref{fig:thr}).
\item
$K=\{1,\ldots,m\}$, where $m=1200$.
\end{itemize}
The solution space $\HH$ is the standard Euclidean space
$\mathbb{R}^{N}$, and \eqref{e:94}
gives rise to the special case of \eqref{e:2}
\begin{equation}
\label{e:sig2}
\text{find}\;\;x\in\HH\;\;\text{such that}\;\;(\forall k\in K)
\quad r_k=Q_{\gamma}\scal{x}{e_k},
\end{equation}
in which $J=\emp$.
\noindent
Note that the standard soft-thresholder on $[-\gamma,\gamma]$ 
can be written as
\begin{equation}
\label{e:st-r}
\soft{\gamma}\colon\xi\mapsto
\sign(Q_{\gamma}\xi)
\left(\sqrt{(Q_{\gamma}\xi)^2+\gamma^2}-\gamma\right).
\end{equation}
To formulate \eqref{e:22} we set, for every $k\in K$, 
\begin{equation}
\label{e:ex2-p}
S_k\colon\RR\to\HH\colon\xi\mapsto
\sign(\xi)\left(\sqrt{\xi^2+\gamma^2}-\gamma\right)e_k,
\end{equation}
which fulfills Assumption~\ref{a:1}\ref{a:1iii} and yields
$S_k\circ R_k=(\soft{\gamma}\scal{\cdot}{e_k})e_k$ \cite{Ibap20}.
We apply \eqref{e:alg} with \eqref{e:lsteps} and the following
control scheme. We split $K$ into $12$ blocks of $100$ consecutive
indices, and select $I_n$ by periodically sweeping through the
blocks, hence satisfying \eqref{e:33} with $M=12$. Moreover,
$\omega_{i,n}\equiv 1/100$. The reconstructed signal shown in
Fig.~\ref{fig:sig} illustrates the ability of the proposed approach
to effectively exploit nonlinearly generated data.

\subsection{Image recovery}
\label{sec:33}
\begin{figure}[t]
\centering
\begin{tabular}{@{}c@{}c@{}}
\includegraphics[width=4.2cm]{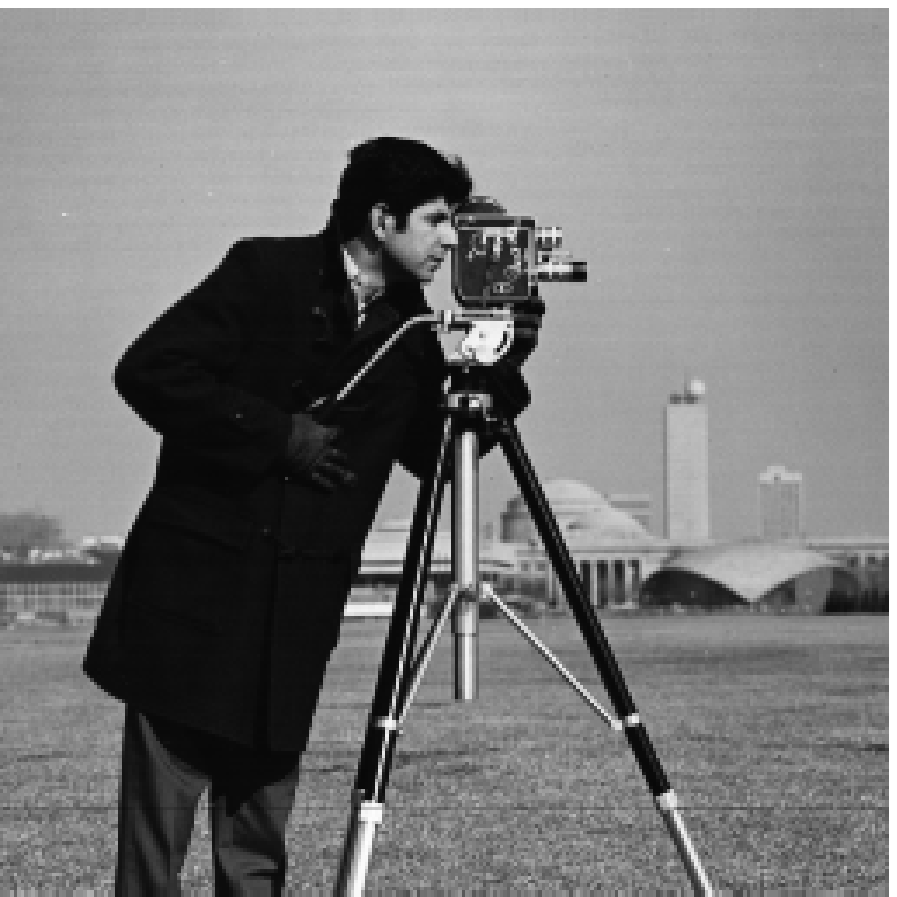}&
\hspace{0.2cm}
\includegraphics[width=4.2cm]{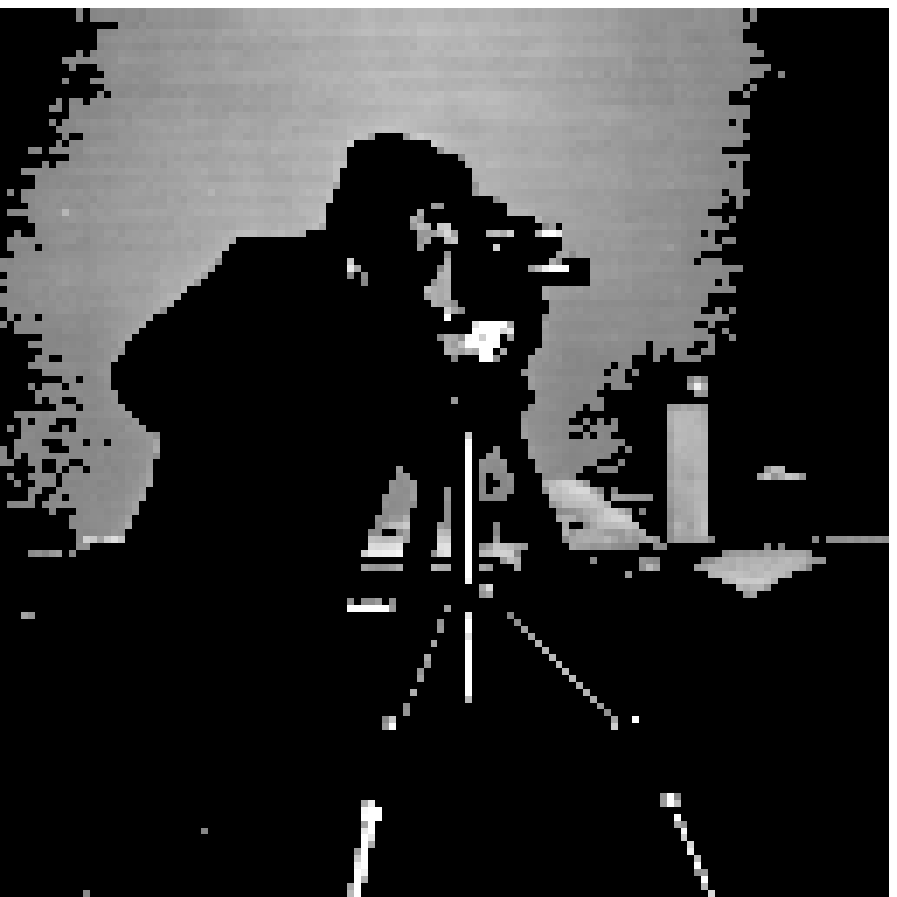}\\
\small{(a)} & \small{(b)}\\
\includegraphics[width=4.2cm]{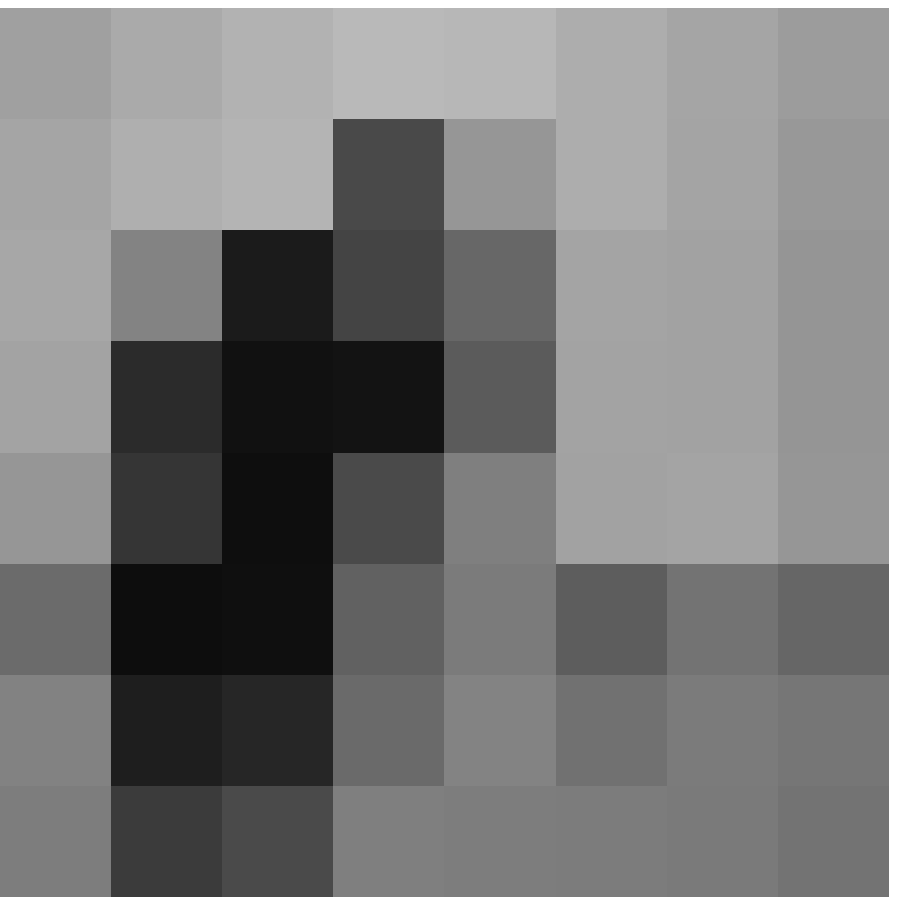}&
\hspace{0.2cm}
\includegraphics[width=4.2cm]{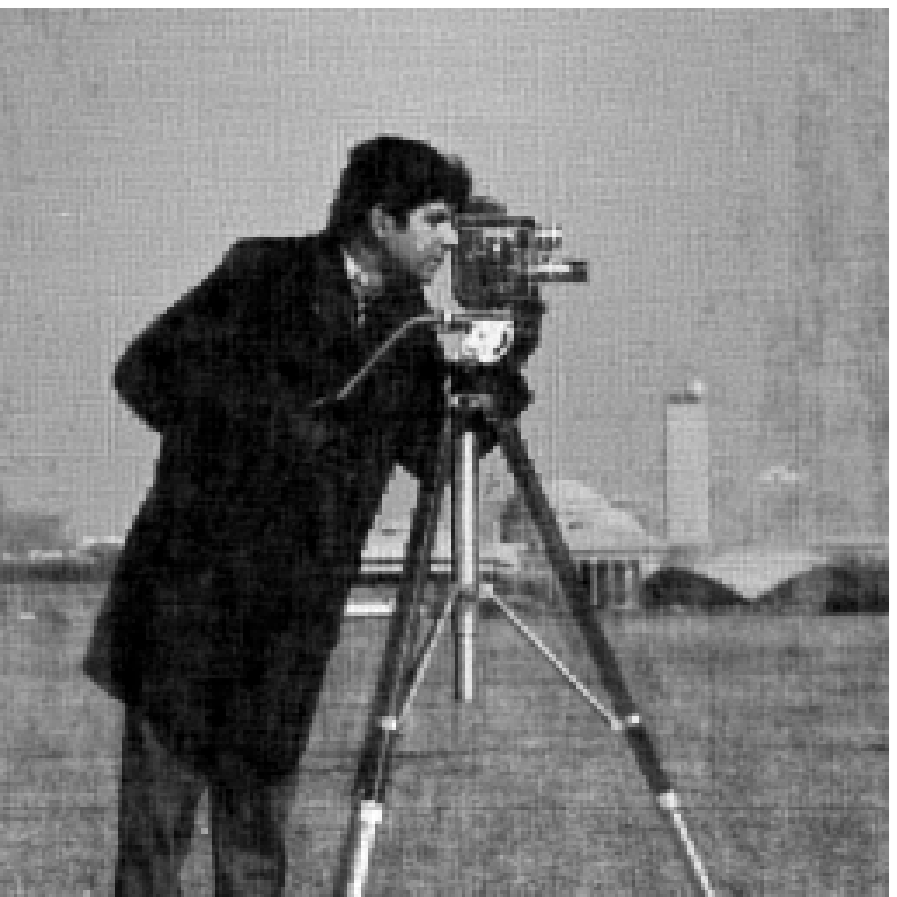}\\
\small{(c)} & \small{(d)}\\
\end{tabular} 
\caption{Images from Section~\ref{sec:33}.
(a) Original image $\overline{x}$.
(b) Compressed image $W^*r_4$.
(c) Down-sampled $8\times 8$ image $r_5$.
(d) Recovered image.}
\label{fig:IM1}
\end{figure}

The goal is to recover the $N\times N$ ($N=256$) image
$\overline{x}$ from the following (see Fig.~\ref{fig:IM1}):
\begin{itemize}
\item 
The Fourier phase $\angle\dft(\overline{x})$
($\dft(\overline{x})$ denotes the 2D Discrete Fourier Transform
of $\overline{x}$).
\item 
The pixel values of $\overline{x}$ reside in $[0,255]$.
\item 
An upper bound $\gamma_3$ on the total variation
$\operatorname{tv}(\overline{x})$ \cite{Imag04}. In this
experiment, $\gamma_3=1.2\operatorname{tv}
(\overline{x})=1.10\times 10^{6}$.
\item 
A compressed representation $r_4=R_4\overline{x}$. Here,
$R_4=Q_4\circ W$, where $W$ is the 2D Haar wavelet transform and
$Q_4$ performs componentwise hard-thresholding via ($\rho=325$)
\begin{equation}
(\forall\xi\in\RR)\quad\hard{\rho}\xi=
\begin{cases}
\xi,&\text{if}\;\;|\xi|>\rho;\\
0,&\text{if}\;\;|\xi|\leq\rho.
\end{cases}
\end{equation}
\item
A down-sampled blurred image $r_5=R_5\overline{x}$. 
Here $R_5=Q_5\circ H_5$, where the linear operator 
$H_5\colon\RR^{N\times N}\to\RR^{N\times N}$ 
convolves with a $5\times 5$ Gaussian kernel with variance $1$,
and $Q_5\colon\RR^{N\times N}\to\RR^{8\times 8}$ maps the average
of each of the 64 disjoint $32\times 32$ blocks of an $N\times N$
image to a representative pixel in an $8\times 8$ image
\cite{Nasr14}. 
\end{itemize}
The solution space is $\HH=\RR^{N\times N}$ equipped with the
Frobenius norm $\|\cdot\|$. To cast the recovery task as an
instance of \eqref{e:2}, we set $J=\{1,2,3\}$, $J_1=\{1,2\}$,
$K=\{4,5\}$,
$C_1=\menge{x\in\HH}{\angle\dft(x)=\angle\dft(\overline{x})}$,
$C_2=[0,255]^{N\times N}$, $f_3=\operatorname{tv}-\gamma_3$,
and $C_3=\menge{x\in\HH}{f_3(x)\leq 0}$. Expressions for
$\proj_{C_1}$ and $\partial f_3$ are provided in \cite{Levi83} and
\cite{Imag04}, respectively. The objective is to
\begin{equation}
\label{e:im3}
\text{find}\;\;x\in\bigcap_{j=1}^3C_j\;\;\text{such that}\;\;
\begin{cases}
R_4x=r_4;\\
R_5x=r_5.
\end{cases}
\end{equation}
Let us verify that Assumption~\ref{a:1}\ref{a:1iii} holds. 
For every $\xi\in\RR$,
\begin{equation}
\label{e:ht-st}
\quad\soft{\rho}\xi =
\hard{\rho}\xi+\begin{cases}
-\rho,&\text{if}\;\;\hard{\rho}\xi>\rho;\\
0,&\text{if}\;\;-\rho\leq\hard{\rho}\xi\leq\rho;\\
\rho,&\text{if}\;\;\hard{\rho}\xi<-\rho.
\end{cases}
\end{equation}

\noindent We construct $S_4$ such that 
$S_4\circ R_4=W^{-1}\circ T\circ W$,
where $T$ applies $\soft{\rho}$ componentwise. In turn, recalling
that $r_4$ is the result of hard-thresholding, $S_4r_4$ is
built by first adding the quantity on the right-hand side of
\eqref{e:ht-st} to $r_4$ componentwise, and then applying the
inverse Haar transform. This guarantees that $S_4$ satisfies
Assumption~\ref{a:1}\ref{a:1iii} \cite{Ibap20}.
Next, we let $D_5\subset\HH$ be the subspace of 
$32\times 32$-block-constant matrices and construct an 
operator $S_5$ satisfying Assumption~\ref{a:1}\ref{a:1iii}
and the identity
$S_5\circ R_5=H_5\circ\proj_{D_5}\circ H_5$ \cite{Ibap20}. In turn,
$S_5r_5=H_5s_5$, where $s_5\in D_5$ is built by repeating
each pixel value of $r_5$ in the block it represents. We thus
arrive at an instance of \eqref{e:22}, which we solve 
using \eqref{e:alg} with \eqref{e:lsteps} and 
\begin{equation}
(\forall n\in\NN)\;\;I_n=J\cup K\;\text{and}\;
(\forall i\in I_n)\;\omega_{i,n}=1/5.
\end{equation}
The resulting image displayed in Fig.~\ref{fig:IM1}(d) shows that
our framework makes it possible to exploit the information from 
the three prior constraints and from the transformations $r_4$ and 
$r_5$ to obtain a quality recovery.

\section{Inconsistent problems}
\label{sec:4}
Inaccuracies and unmodeled dynamics may cause \eqref{e:2} to admit
no solution. In such instances, we propose the following
relaxation for \eqref{e:2} \cite{Ibap20}.

\begin{assumption}
\label{a:r}
For every $j\in J$, the operator $\proj_{C_j}$ is
easily implementable and, for every $k\in K$,
Assumption~\ref{a:1}\ref{a:1iii} holds. In addition, 
$\{\omega_j\}_{j\in J}\subset\rzeroun$ and 
$\{\omega_k\}_{k\in K}\subset\rzeroun$ satisfy 
$\sum_{j\in J}\omega_j+\sum_{k\in K}\omega_k=1$. 
\end{assumption}

\noindent Under Assumption~\ref{a:r}, the goal is to
\begin{multline}
\label{e:r}
\text{find}\;\;x\in\HH\;\;\text{such that}\\
\sum_{j\in J}\omega_j(x-\proj_{C_j}x)+\sum_{k\in K}
\omega_k(S_kR_kx-S_kr_k)=0.
\end{multline}
When $K=\varnothing$, the solutions of \eqref{e:r} are the 
minimizers of the least squared-distance proximity function
$\sum_{j\in J}\omega_j d^2_{C_j}$ \cite{Aiep96}. If
\eqref{e:2} does have solutions, then it is equivalent to
\eqref{e:r}. The algorithm of \cite{Comb20} can be used to solve
\eqref{e:r} block-iteratively.

\end{document}